\documentclass[journal]{IEEEtran}

\usepackage{mathtools}
\usepackage{bm}
\usepackage{algorithmicx}
\usepackage{algorithm,setspace}
\usepackage{algpseudocode}
\usepackage{graphicx}
\usepackage{epstopdf}
\usepackage{amsfonts}
\usepackage{amsmath}
\usepackage{multirow}
\usepackage{float}
\usepackage{epsfig}
\usepackage{xcolor}
\usepackage{amssymb}
\usepackage[caption=false]{subfig}

\newtheorem{theorem}{Theorem}[section]

\newtheorem{lemma}[theorem]{Lemma}

\newtheorem{assumption}[theorem]{Assumption}

\newtheorem{remark}[theorem]{Remark}

\usepackage{balance}

\usepackage{colortbl}
\newcommand\oprocendsymbol{\hbox{$\square$}}
\newcommand\oprocend{\relax\ifmmode\else\unskip\hfill\fi\oprocendsymbol
}

\usepackage{tikz}
\usepackage{url}
\definecolor{blue@O4S}{RGB}{0, 41, 69}
\definecolor{emph@O4S}{RGB}{0, 93, 137}
\definecolor{red@O4S}{RGB}{127,0,0}
\definecolor{gray@O4S}{RGB}{112, 112, 112}

\usepackage[printonlyused]{acronym}

\newcommand{\subj}{\text{subj. to}}

\newcommand{\argmin}{\mathop{\rm argmin}}
\newcommand{\argmax}{\mathop{\rm argmax}}

\newcommand{\real}{\mathbb{R}}
\newcommand{\until}[1]{\{1,\ldots,#1\}} 

\renewcommand{\natural}{\mathbb{N}}

\newcommand{\Nag}{N}

\newcommand{\EE}{\mathcal{E}}
\newcommand{\GG}{\mathcal{G}}

\newcommand{\innbrs}{\mathcal{N}^{\textrm{\textnormal{in}}}}

\newcommand{\DASS}{Distributed Dynamic Assignment and Servicing Strategy}

\newcommand{\ipmp}{IP-MP} %

\newcommand{\1}[1]{1_{#1}}
\newcommand{\0}[1]{0_{#1}}

\newcommand{\varz}{z}

\newcommand{\indexa}{\ell}
\newcommand{\indexb}{m}

\newcommand{\PP}{P}

\algnewcommand\algorithmicinit{\textbf{Initialization:}}
\algnewcommand\init{\item[\algorithmicinit]}
\algnewcommand\algorithmicawake{\textsf{\textit{{AWAKE}}}}
\algnewcommand\awake{\item[\algorithmicawake]}
\algnewcommand\algorithmicidle{\textsf{\textit{{IDLE}}}}
\algnewcommand\idle{\item[\algorithmicidle]}
\algnewcommand\algorithmicinput{\textbf{Input:}}
\algnewcommand\inpt{\item[\algorithmicinput]}
\algnewcommand\algorithmicoutput{\textbf{Output:}}
\algnewcommand\outpt{\item[\algorithmicoutput]}

\newcommand{\EEt}{\EE^t}
\newcommand{\GGt}{\GG^t}
\newcommand{\wel}{w}
\newcommand{\pel}{p}
\newcommand{\vel}{g}
\newcommand{\Ntask}{M}
\newcommand{\varx}{x}

\newcommand{\zinc}[2]{z_{{[#1]}}^{{#2}}}
\newcommand{\extr}{v^q}
\newcommand{\mult}{\lambda^q}
\newcommand{\QQ}{Q}
\newcommand{\barQQ}{\bar{Q}}

\newcommand{\cc}{c}
\newcommand{\DD}{D}

\newcommand{\basisarg}[2]{B_{{#1}}^{{#2}}}
\newcommand{\colg}{h_i}
\newcommand{\Prob}{\mathcal{MP}}

\newcommand{\jarg}[2]{J_{{#1}}^{{#2}}}
\newcommand{\costval}{J}
\newcommand{\unif}[2]{\mathcal{U}(#1,#2)}
\newcommand{\labelarg}[2]{\mathcal{L}_{#1}^{#2}}
\newcommand{\convT}{\overline{T}}
\newcommand{\lambdaRMP}{\tilde{\Lambda}}
\newcommand{\lambdavarRMP}{\tilde{\lambda}}
\newcommand{\lambdaell}{\Lambda^{\!\star\ell}}
\newcommand{\tree}{\mathcal{T}}
\newcommand{\treearg}[2]{\tree_{{#1}}^{{#2}}}
\newcommand{\varzLP}{\varz_{[i]}^{\textsc{LP}}}
\newcommand{\JLP}{J_{i}^{\textsc{LP}}}
\newcommand{\solfunc}{\textsc{ExtractSol}}
\newcommand{\branchfunc}{\textsc{Branch}}
\newcommand{\upfunc}{\textsc{ExtractConstr}}
\newcommand{\varzell}{\varz^{\star\ell}}
\newcommand{\jell}{J^{\star\ell}}
\newcommand{\bart}[1]{\bar{t}_{#1}}
\newcommand{\probarg}[2]{\Prob_{#1}^{#2}}

\acrodef{GAP}{\emph{Generalized Assignment Problem}} %
\acrodef{RMP}{\emph{Restricted Master Problem}} %

\graphicspath{{figs/}}

\makeatletter
\newcommand{\StatexIndent}[1][3]{%
  \setlength\@tempdima{\algorithmicindent}%
  \Statex\hskip\dimexpr#1\@tempdima\relax}
\makeatother

\usepackage{cuted}
\begin{document}

\title{Generalized Assignment for Multi-Robot Systems via Distributed Branch-And-Price}
\author{Andrea Testa and %
Giuseppe Notarstefano \thanks{This result is part of a project that has received funding from the
European Research Council (ERC)
    under the European Union's Horizon 2020 research and innovation programme
    (grant agreement No 638992 - OPT4SMART). }
\thanks{A. Testa and G. Notarstefano are  with the Department of Electrical, Electronic and Information Engineering,
University of Bologna, Bologna, Italy, \texttt{a.testa@unibo.it, giuseppe.notarstefano@unibo.it}.}
}

\maketitle
\begin{strip}\leavevmode\kern15pt
\begin{minipage}{\dimexpr\linewidth-30pt\relax}
{\vspace{-2.5cm}
\bf \textcopyright 2021 IEEE. Personal use of this material is permitted.  Permission from IEEE must be obtained for all other uses, in any current or future media, including reprinting/republishing this material for advertising or promotional purposes, creating new collective works, for resale or redistribution to servers or lists, or reuse of any copyrighted component of this work in other works.}
\end{minipage}
\end{strip}
\begin{abstract}
In this paper, we consider a network of agents that has to self-assign a set of tasks while respecting resource constraints. One possible formulation is the Generalized Assignment Problem, where the goal is to find a maximum payoff while satisfying capability constraints. We propose a purely distributed branch-and-price algorithm to solve this problem in a cooperative fashion. Inspired by classical (centralized) branch-and-price schemes, in the proposed algorithm each agent locally solves small linear programs, generates columns by solving simple knapsack problems, and communicates to its neighbors a fixed number of basic columns. 
We prove finite-time convergence of the algorithm to an optimal solution of
  the problem.  Then, we apply the proposed scheme to a generalized assignment
  scenario in which a team of robots has to serve a set of tasks.
We implement the proposed algorithm in a ROS testbed and provide experiments for a team of heterogeneous robots solving the assignment problem.
\end{abstract}

\IEEEpeerreviewmaketitle

\section{Introduction}
The \ac{GAP} is a well known combinatorial optimization problem with several applications as
vehicle routing, facility location,
resource scheduling and supply chain, to name a few~\cite{oncan2007survey,savelsbergh1997branch,fisher1981generalized}.
Even though \ac{GAP} is a NP-hard problem, several approaches
have been developed for solving this problem both for exact and approximate solutions.
We refer the reader to~\cite{martello1992generalized} for a survey.
Branch-and-price algorithms~\cite{savelsbergh1997branch,barnhart1998branch}
are among the most investigated algorithms allowing for both optimal
and suboptimal solutions.

\subsection{Related Work}

Task assignment naturally arises in cooperative robotics, where heterogeneous agents collaborate
to fulfill a complex task, see, e.g.,~\cite{gerkey2004formal} for an early reference.
Specific applications include persistent monitoring of locations~\cite{hartuv2018scheduling}, path planning of mobile robots, e.g., UAVs~\cite{bellingham2003multi}, task scheduling for robots working in the same space~\cite{gombolay2018fast}, vehicle routing~\cite{turpin2015approximation} and
task assignment in urban environments~\cite{shima2007assigning}.
All the previous problems are solved by means of centralized approaches.

In order to deal with the computational complexity of the problem, a branch of literature
analyzes parallel and decentralized approaches to the problem\footnote{We denote \emph{parallel} the
approaches based on master-slave architectures, while we call \emph{decentralized} the schemes with
independent agents that do not communicate among each other.}.
A well known parallel approach is the
auction based one, originally proposed in~\cite{bertsekas1988auction}.
A market-based approach is considered in~\cite{dias2008sliding} for the coordination of human-robot teams.
Authors in~\cite{castanon2003distributed} propose an
algorithm, based on a
sequential shortest augmenting path scheme, to solve a dynamic multi-task
allocation problem. Agents propose assignments
that are validated by a coordinating unit.
As for decentralized schemes,
authors in~\cite{lerman2006analysis} solve a dynamic task allocation problem
for robots that can perform local sensing operations and do not
communicate among each other.
In~\cite{alighanbari2005decentralized}, a task
assignment problem is solved, in a decentralized scheme, through the so called petal algorithm. In~\cite{nam2019robots}, a dynamic task assignment
problem in which the cost vector changes in a bounded region is considered. A central unit is initially required but robots are able to exploit local communications to perform a reallocation if needed.
An area partitioning problem for multi-robot systems is proposed in~\cite{hassan2014task} and solved by a genetic algorithm. An area coverage problem in marine environments is solved in~\cite{karapetyan2018multi} with heuristics based on the traveling salesman problem.

As for distributed schemes, i.e., with processors in a peer-to-peer network
  without a central coordinator, a distributed version of the Hungarian method
  is proposed in~\cite{chopra2017distributed}.  A distributed simplex scheme
for degenerate linear programs (LP) is proposed in~\cite{burger2012distributed}
in the context of multi-agent assignment problems.
A distributed subgradient is applied in~\cite{settimi2013subgradient}
to a task assignment problem, while in~\cite{burger2011locally} a distributed column generation scheme is proposed.
A linear task assignment
problem with time-varying cost functions is considered in~\cite{montijano2019distributed}, %
while an optimal role and position assignment problem is addressed in \cite{mosteo2017optimal}
by iteratively solving a sequence of linear assignment problems.
In these approaches, authors
neglect integrality constraints on the decision variables, relying on the unimodular structure of the problems.
As for distributed, suboptimal approaches for task assignment problems,
in~\cite{karaman2008large} a large-scale distributed task/target assignment
problem across a fleet of autonomous UAVs is considered, but the
  communication graph is assumed to be complete.
In the context of wireless sensor networks,~\cite{pilloni2016deployment}
proposes
a distributed task allocation in order to maximize the network life-time.
A distributed task assignment algorithm is used in conjunction with a
deterministic annealing in~\cite{kwok2011distributed}, in the context of
limited-range sensor coverage.  In~\cite{abbatecola2018distributed}, a dynamic
vehicle routing problem is approached with a distributed protocol in which
agents iteratively solve graph partitioning
  problems. %
In the works~\cite{testa2017finite,testa2019distributed} authors address
  Mixed-Integer Linear Programs by means of a distributed cutting-plane
  algorithm and apply it to a multi-agent multi-task assignment problem.
Distributed implementations of the auction-based algorithm are often used to solve task assigment problems,
see, e.g.,~\cite{choi2009consensus} for an early reference, and in particular GAPs~\cite{luo2013distributed,luo2015distributed}. The auction-based approach allows for a suboptimal solution with performance guarantees.
In~\cite{williams2017decentralized}, this approach is applied to a task
allocation problem expressed as a combinatorial optimization problem with
matroid constraints.  In the recent
works~\cite{buckman2019partial,talebpour2019adaptive}, a dynamic task allocation
scenario with partial replanning is considered.
The above references show that the exact resolution of
\ac{GAP} is an open problem in a purely distributed setting.
Indeed, state-of-the-art solutions are usually based on proper linear
relaxations or suboptimal approaches.
\subsection{Contributions}
In this paper, we propose a purely distributed version of the
branch-and-price algorithm to solve the Generalized Assignment Problem by means
of a network of agents.
Specifically, each agent locally solves a linear programming relaxation of the
\ac{GAP}, generates columns by solving a (simple) knapsack problem, and
exchanges estimates of the solution with neighboring agents.  Due to the
relaxation of the integrality constraints, the solution of this problem may not
be feasible. Thus, new problems, based on the original one with suitable
additional constraints, have to be solved. The set of these problems can be
represented by a so called branching tree. By leveraging on their communication
capabilities, agents explore their local trees until an optimal solution of the
optimization problem has been found.
With respect to the aforementioned works, and
  specifically~\cite{karaman2008large,choi2009consensus,luo2013distributed}, the
  proposed scheme has the following distinctive new features. To the best of the
  authors' knowledge, this is the first attempt to solve \ac{GAP} to optimality
  in a purely distributed fashion.  Remarkably, the proposed scheme is shown to
  converge also under time-varying and directed communication
  networks. Moreover, it is worth noticing that the approaches
  in~\cite{testa2017finite,testa2019distributed} are not applicable to the
  considered GAP scenario, which involves equality constraints.
  Finally, we apply the proposed algorithm to a dynamic task assignment problem
  where tasks may arrive during time and robots have to adapt the local plan
  according to the new information.
An experimental platform, based on ROS (Robot Operating System), is proposed to
run experiments in which a team of aerial and ground robots cooperatively
solve the \ac{GAP} relying on the proposed distributed branch-and-price scheme.

The paper unfolds as follows.
In Section~\ref{sec:CentrGAP} we introduce the distributed setup considered throughout the paper. Then, we introduce the Generalized Assignment Problem
and a centralized scheme, called \emph{branch-and-price}, to solve it.
In Section~\ref{sec:Distributed} we propose a purely distributed branch-and-price algorithm. %
In Section~\ref{sec:Numerical} we provide numerical simulations for randomly generated GAPs and in Section~\ref{sec:experiment} we show the results of experiments on a swarm of heterogeneous robots.%

\paragraph*{Notation}
We denote by $e_\indexa$ the $\indexa$-th vector of the canonical basis (e.g.,
$e_1 = [1\, 0\, \ldots\, 0]^\top$) of proper dimension.
Given a vector $v_\indexa \in \real^d$, we denote by $v_{\indexa_\indexb}$ the $\indexb$-th
component of $v_\indexa$. Also, we denote by $\1{r}$ ($\0{r}$) the vector in $\real^{r}$ with all
its entries equal to $1$ ($0$).
\section{Distributed Setup and Preliminaries}

\label{sec:CentrGAP}
In this section, we introduce the distributed setup for the Generalized Assignment
Problem addressed in the paper. Also,
the (centralized) branch-and-price scheme is illustrated.

\subsection{Distributed Problem Setup}

In the Generalized Assignment Problem,
the objective is to find a maximal profit assignment of $\Ntask$ tasks to
$\Nag$ agents such that
each task is assigned only to one agent.
In this scenario, the generic agent $i$ has a reward
$\pel_{i\indexb}\in\real$ if it executes the $\indexb$-th task.
It also has a limited capacity $\vel_i\in\real$ and it uses
an amount $\wel_{i\indexb}\in\real$ of capacity if it performs
the $\indexb$-th task. Let $\varx_{i\indexb}$ be a binary variable indicating
whether task $\indexb$ is assigned to agent $i$ ($x_{i\indexb} =
1$) or not ($x_{i\indexb} = 0$). We denote constraints in the form
$\varx_{i\indexb}\in\{0,1\}$ as integer
constraint.
Then, the standard integer programming formulation is %
\begin{align} \label{eq:GAP}
\begin{split}
\max_{x_{11},\ldots,x_{\Nag\!\Ntask}} \: & \: \sum_{i =1}^\Nag \sum_{\indexb =1}^\Ntask
\pel_{i\indexb} \varx_{i\indexb}
  \\
\subj \: & \: \sum_{i=1}^\Nag \varx_{i\indexb} = 1, \indexb= 1,\ldots,\Ntask,
  \\
& \: \sum_{\indexb=1}^\Ntask \wel_{i\indexb} x_{i\indexb} \leq \vel_i,  i= 1,\ldots,\Nag,
  \\
& \: x_{i\indexb} \in \{0,1\}, i= 1,\ldots,\Nag,  \indexb= 1,\ldots,\Ntask.
\end{split}
\end{align}
In order to streamline the notation, we now introduce
a formulation of the \ac{GAP} better highlighting the structure
of the problem in a distributed scenario.
Let $\varz_i =
[\varx_{i1},\ldots,\varx_{i\Ntask}]^\top \in \real^\Ntask, \forall i=1,\ldots,\Nag$.
In the following we denote as $z$ the stack $[\varz_1^\top,\ldots,\varz_{\Nag}^\top]^\top$.
Also,
let $\cc_i= [\pel_{i1},\ldots,\pel_{i\Ntask}]^\top\in \real^\Ntask $,
$\DD_i= [\wel_{i1},\ldots,\wel_{i\Ntask}]\in \real^\Ntask $ and $\PP_i = \{ \varz_i \in \{0, 1\}^{\Ntask} \mid \DD_i  \varz_i \leq \vel_i\}$,
for $i = 1, \ldots, \Nag$.
Then, \eqref{eq:GAP} can be recast as
\begin{align} \label{eq:GAP2}
\begin{split}
  \max_{\varz_1,\ldots,\varz_\Nag} \: & \: \sum_{i =1}^\Nag c_i^\top \varz_i
  \\
  \subj \: & \: \sum_{i=1}^\Nag \varz_i = \1{\Ntask},
  \\
  & \: \varz_i \in \PP_i, i= 1,\ldots, \Nag.
  \\
\end{split}
\end{align}
This new formulation of~\eqref{eq:GAP} allows us to point out the distributed nature of the problem. Namely, $c_i$ describes the profits
associated to assigning tasks to agent $i$,
$\sum_{i=1}^\Nag \varz_i = \1{\Ntask}$ describes the assignment constraints
(\emph{coupling constraints}),
$\PP_i$
describes the capacity restrictions on the agents (\emph{local constraints}).
It is worth noting that the sets $\PP_i$ are bounded
for $i=1,\ldots, \Nag$.

The agents must solve~\eqref{eq:GAP2} cooperatively in a distributed fashion
with limited
communication and computation capabilities, as well as limited memory.
We consider the natural scenario in which the $i$-th agent only
knows the polyhedron $P_i$
and the cost vector
$c_i$, thus not having knowledge of other agent data.
In order to solve the problem, agents can exchange information according to a
time-varying communication network modeled as a time-varying digraph
$\GG^t=(\until{N},\EE^t)$, with $t\in\natural $ being a universal slotted
time representing a temporal information on the graph evolution.
Notice that time $t$ does not need to be known by the agents.
A digraph $\GGt$ models the communication in the
sense that there is an edge $(i,j) \in \EEt$ if and only if agent $i$ is able
to send information to agent $j$ at time $t$.
For each node $i$, the set of
\emph{in-neighbors} of $i$ at time $t$ is denoted by $\innbrs_{i,t}$ and is the
set of $j$ such that there exists an edge $(j,i) \in \EEt$.
A static digraph is said to be \emph{strongly connected} if there exists a
directed path for each pair of agents $i$ and $j$.
Next, we require the following.
\begin{assumption}[Graph Connectivity]\label{ass:connectivity}
The communication graph is $L$-strongly connected, i.e.,
there exists an integer $L\geq1$ such that, for all $t\in\natural$, the graph
 $(\until{\Nag},\bigcup_{\tau=t}^{t+L-1}\EE^{\tau})$ is strongly connected.\oprocend
\end{assumption}
Notice that this is a standard, mild, assumption in the context of
distributed optimization that allows to model direct,
time-varying, asynchronous and possibly unreliable communication.

\subsection{Centralized Branch-and-Price Method}
We now introduce the main concepts regarding the
branch-and-price scheme.
We refer the reader to~\cite{barnhart1998branch}
for a more detailed dissertation.
For the sake of clarity, we organize this subsection in three parts.

\subsection*{Dantzig-Wolfe Decomposition for GAPs}
An equivalent formulation of~\eqref{eq:GAP2}, which is exploited in our distributed setup, can be
obtained as follows.
Such procedure, originally introduced in~\cite{barnhart1998branch},
is strictly related to the
\emph{Dantzig-Wolfe Decomposition}~\cite{vanderbeck2000dantzig}.
Points $\varz_i\in\PP_i$ can be represented as the linear combination of a finite number of vectors $\extr_i$, with $q\in\until{|\QQ_i|}$, i.e.,
\begin{align}
\label{eq:comb_zi}
\begin{split}
  \varz_i = &\sum_{q=1}^{|\QQ_i|} \extr_i \mult_i,
\end{split}
\end{align}
where we denote with $\QQ_i$ the set of vectors $\extr_i$.
The variables $\mult_i$, $\forall q= 1,\ldots,|\QQ_i|$, also called \emph{combiners}, have to satisfy
\begin{align}
\label{eq:comb_con}
\begin{split}
 &\sum_{q=1}^{|\QQ_i|} \mult_i = 1,
  \\
  &\mult_i \in \{0, 1\}. %
\end{split}
\end{align}
It can be shown that, for \ac{GAP}s, $\QQ_i$ coincides with
the set of extreme points of the convex hull $\text{conv}(\PP_i)$ of $\PP_i$,~\cite{barnhart1998branch}.
Let $\Lambda\in\real^{\sum_{i=1}^N|\QQ_i|}$, be
the stack of all the combiners.
Substituting~\eqref{eq:comb_zi} and~\eqref{eq:comb_con} in~\eqref{eq:GAP2} leads to the following equivalent Integer Programming Master Problem (\ipmp{})
\begin{align} \label{eq:GAP_IPMP}
\begin{split}
\max_{\Lambda} \: & \: \sum_{i =1}^\Nag \sum_{q=1}^{|\QQ_i|} (\cc_i^\top \extr_i)\mult_i
  \\
\subj \: & \: \sum_{i=1}^\Nag \sum_{q=1}^{|\QQ_i|} \extr_i \mult_i =\1{\Ntask},
  \\
  & \: \sum_{q=1}^{|\QQ_i|} \mult_i = 1,  i= 1,\ldots, \Nag,
  \\
  & \: \mult_i \in \{0,1\},  q\in\until{|\QQ_i|}, i= 1,\ldots, \Nag.
\end{split}
\end{align}

Notice that an optimal solution $\varz^\star$ of~\eqref{eq:GAP2} can be retrieved from an optimal solution $\Lambda^{\star}$ of~\eqref{eq:GAP_IPMP} by substituting the entries of $\Lambda^{\!\star}$ in~\eqref{eq:comb_zi},~\cite{savelsbergh1997branch}.

\subsection*{Branching Tree}
\begin{figure}
	\centering
	\includegraphics[scale=2]{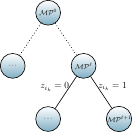}
	\caption{Example of branching tree with problems generated according to constraints on the sets $\PP_i$.}
	\label{fig:tree}
\end{figure}

The presence of binary constraints makes~\eqref{eq:GAP_IPMP} hard to solve. In order to
find an optimal solution to the problem, the branch-and-price algorithm,~\cite{barnhart1998branch}, explores the set of feasible solutions of GAP by iteratively generating and solving
\emph{relaxed} versions of~\eqref{eq:GAP_IPMP} including suitable, tightening constraints.
That is, the constraint $\mult_i \in \{0,1\}$ of all these problems is relaxed to
$\mult_i \geq 0$ ($\mult_i \leq 1$ can be omitted as
it is implicit in the constraint $\sum_{q=1}^{|\QQ_i|} \mult_i=1$).
These problems can be represented as nodes of a so called \emph{branching tree}, see, e.g., Figure~\ref{fig:tree}.
The $\ell$-th node of the tree to be solved represents a problem, in the form of~\eqref{eq:GAP_IPMP}, obtained by relaxing the integer constraints and enforcing the constraints of the edges. We denote such problem as $\Prob^\ell$, by $\lambdaell$ its optimal solution
and by $\varzell$ and
$\jell$ the solution and cost in terms of the variables $\varz$.
Branches (edges) indicate the constraints that have to be added to generate new problems. %
Instead of considering constraints on the $\Lambda$ variables, new problems are generated by including, to the sets $\PP_i$, additional constraints in the form $\varz_{i_k}=0, \varz_{i_k}=1$ for some $\varz_{i_k}(\ell)\not\in\{0,1\}$.
In the following, we assume there exists an extraction strategy to determine the next node of the tree to be solved and a rule to choose the index $i_k$.
Let $\PP_i^\ell$ be the sets obtained by including such additional constraints and let $\QQ_i^\ell$ be the sets of extreme points of $\text{conv}(\PP_i^\ell)$. Since $\QQ_i^\ell\subset\QQ_i$, this results in generating relaxed problems in the form of~\eqref{eq:GAP_IPMP} with less optimization variables.
The algorithm keeps track of a \emph{lowerbound}, also called \emph{incumbent}, $J^{\textsc{inc}}$
on the cost, and of a candidate solution $\varz^{\textsc{inc}}$.
After solving the generic problem $\Prob^\ell$, one of the following operations is performed:
\begin{itemize}
\item \emph{Incumbent Update:}
If $\varzell\in \{0,1\}^{\Nag\Ntask}$ and $\jell\geq J^{\textsc{inc}}$, then
$J^{\textsc{inc}}=\jell$ and $\varz^{\textsc{inc}}=\varzell$.
\item \emph{Branching:} If $\jell > J^{\textsc{inc}}$ and $\varzell\not\in \{0,1\}^{\Nag\Ntask}$
and $\jell\geq J^{\textsc{inc}}$, two new problems are
added to the tree.
\item \emph{Pruning:} If $\jell\leq J^{\textsc{inc}}$ or $\Prob^\ell$ is infeasible, nothing is done.
\end{itemize}
\begin{remark}
Pruning prevents the algorithm from inspecting
problems that do not improve $J^{\textsc{inc}}$.~\oprocend
\end{remark}
At the end of the algorithm $J^{\textsc{inc}}$ and $\varz^{\textsc{inc}}$ coincide with the
optimal cost $J^\star$ and
solution
$\varz^\star$ of~\eqref{eq:GAP2}.
\subsection*{Column Generation}
Each problem $\Prob^\ell$ has a large number of optimization variables. Thus, it can be approached by means of the so called \emph{column generation algorithm} (originally
proposed in~\cite{gilmore1961linear} in the context of cutting stock problems).
It consists in iteratively performing the
  following three steps:

i) A \ac{RMP} is solved, made by a small subset of the
\emph{columns}\footnote{We refer the reader to
Appendix~\ref{sec:lps} for the
definition of column for a linear program.
Informally, a column is a portion
of the cost and constraint vectors associated to a decision variable.}
of $\Prob^\ell$ with sets $\barQQ_i^\ell \subset Q_i^\ell$
and a smaller combiner vector $\lambdaRMP$, i.e.,
\begin{align} \label{eq:RMP}
\begin{split}
\max_{\lambdaRMP} \: & \: \sum_{i =1}^\Nag \sum_{q=1}^{|\barQQ_i^\ell|}%
(\cc_i^\top \extr_i)\lambdavarRMP_i^q
  \\
\subj \: & \: \sum_{i=1}^\Nag \sum_{q=1}^{|\barQQ_i|}  \extr_i \lambdavarRMP_i^q = \1{\Ntask},
  \\
& \: \sum_{q=1}^{|\barQQ_i^\ell|}  \lambdavarRMP_i^q = 1,  i= 1,\ldots, \Nag,
  \\
  & \: \lambdavarRMP_i^q\geq0,  q\in\until{|\barQQ_i^\ell|}, i= 1,\ldots, \Nag.
\end{split}
\end{align}
\indent ii) If possible, new columns are added to the \ac{RMP} in
  order to improve the current solution.
Let $[\pi^\top \, \mu^\top]^\top$ be a dual optimal solution of \eqref{eq:RMP} at
a generic iteration of the algorithm. In
particular, $\pi\in\real^\Ntask$ is associated to the constraint $
\sum_{i=1}^\Nag \sum_{q \in \barQQ_i^\ell}  \extr_i \lambdavarRMP_i^q\allowbreak =
\1{\Ntask}$, while $\mu\in\real^\Nag$ is associated to the constraint $\sum_{q
\in \barQQ_i^\ell} \lambdavarRMP_i^q\allowbreak = 1$.
For each $i= 1,\ldots, \Nag$ a new column is found by solving the
so called \emph{pricing problem}:
\begin{align} \label{eq:BinarySubProblem}
\begin{split}
 \bar{v}_i \in \argmax_{\varz_i} \: & \: (\cc_i-\pi)^\top \varz_i
  \\
  \subj \: & \: \varz_i \in \PP_i^\ell.
\end{split}
\end{align}
Consider now the associated column in the form $\colg = [c_i^\top  \bar{v}_i \;\;
\bar{v}_i^\top \;\; e_i^\top]^\top$.
Then, $\colg$ allows for a cost improvement if it has positive reduced
cost, i.e., if
$
(\cc_i - \pi ) \bar{v}_i - \mu_i >0
$.

iii) A \emph{pivoting} operation is performed, that is, all columns with
  positive reduced cost are included in the \ac{RMP}, while columns of the
  \ac{RMP} that are not associated to basic variables are dropped\footnote{We refer the reader to Appedix~\ref{sec:lps} for the definition of basic variables.}.
Then, the procedure is iterated
until no more columns with positive reduced cost can be found.
Let $\lambdaell$ be
the final optimal solution
of the
relaxed version of~\eqref{eq:GAP_IPMP} obtained with this procedure.
Let $\bar{\lambda}^q_i$ be the entry of $\lambdaell$ associated to a
vertex $\extr_i \in\QQ_i^\ell$.
Then, the solution $\varzell$
of $\Prob^\ell$ can be expressed as (c.f.~\eqref{eq:comb_zi})
\begin{align}\label{eq:reconstr_varz}
\varz_i(\ell)=\sum_{q=1}^{|\QQ_i^\ell|}\bar{\lambda}_i^q \extr_i,\quad
i=1,\ldots,\Nag.
\end{align}
\begin{remark}\label{rmk:pric_GAP}
When applying the Dantzig-Wolfe Decomposition,
we follow the approach
in~\cite{barnhart1998branch} and do not relax the binary constraints in~\eqref{eq:GAP2}.
This results in local knapsack problems~\eqref{eq:BinarySubProblem} that can be efficiently solved through dynamic programming schemes,~\cite{martello1990knapsack}.~\oprocend
\end{remark}
We collect in Table~\ref{tb:symbols} all the relevant symbols.
\begin{table}[t]\centering
	\caption{List of the main symbols and their definitions}
	\label{tb:symbols}
	\begin{tabular}{ll}
	  \hline
		\multicolumn{2}{c}{Standard GAP Formulation}
		\\
		\hline
		$\Nag\in\natural_{> 0}$ 
		& Number of robots
		\\
		$\Ntask\in\natural_{> 0}$
		& Number of tasks
		\\
		$x_{im}\in\{0,1\}$
		& $1$ if robot $i$ serves task, $0$ otherwise
		\\
		$\varz_i\in \real^\Ntask$
		& $[\varx_{i1},\ldots,\varx_{i\Ntask}]^\top$
		\\
		$p_{im} \in \mathbb{R}_{\geq0}$
		& Reward if robot $i$ serves task $m$
		\\
		$c_{i} \in \real^M$
		& $[\pel_{i1},\ldots,\pel_{i\Ntask}]^\top$
		\\
		$g_i\in\real_{\geq0}$
		& Capacity of robot $i$
		\\
		$w_{im}\in \mathbb{R}_{\geq0}$
		& Capacity consumption of task $m$ for robot $i$
		\\
		$D_{i} \in \real^M$
		& $[\wel_{i1},\ldots,\wel_{i\Ntask}]$
		\\
		$\PP_{i} \subseteq \{0,1\}^M$
		& $\{ \varz_i \in \{0, 1\}^{\Ntask} \mid \DD_i  \varz_i \leq \vel_i\}$
		\\
    \hline
		\multicolumn{2}{c}{Dantzig-Wolfe Reformulation}
		\\
		\hline
		$Q_i$
		& set of extreme points of $P_i$
		\\
		$\extr_i\in\{0,1\}^M$
		& $q$-th extreme point of $P_i$
		\\
		$\mult_i\in\{0,1\}$
		& Combiner associated to $\extr_i$
		\\
		$\Lambda\in\real^{\sum_{i=1}^N|\QQ_i|}$
		& Stack of combiners
		\\
    	\hline
		\multicolumn{2}{c}{Branch-and-Price}
		\\
		\hline
		$\Prob^\ell$
		& $\ell$-th node of the branching tree
		\\
		$\varzell, \jell$
		& Optimal solution and cost of $\Prob^\ell$
		\\
		$\PP_i^\ell$
		& Constraint set of robot $i$ at node $\Prob^\ell$
		\\
		$\varz^{\textsc{inc}}$, $J^{\textsc{inc}}=\jell$
		& Candidate GAP solution and cost
		\\
           \hline
	\end{tabular}
\end{table}

\section{Distributed Branch-and-Price Method}
\label{sec:Distributed}
In this section, we provide a purely distributed algorithm, inspired by the
  centralized branch-and-price scheme, to solve~\eqref{eq:GAP2} in a
  peer-to-peer network.
We assume a solver for Linear Programs is available.  In particular, we use the
simplex algorithm proposed in \cite{jones2007lexicographic} to find the unique
\emph{lexicographically minimal optimal} solution
of a LP and the associated optimal basis.
In the proposed distributed algorithm, called \emph{Distributed Branch-and-Price},
each agent $i$ maintains and updates, at the generic time $t$, local optimal cost and solution candidates $\jarg{i}{t}$ and $\zinc{i}{t}$, as well as a local tree $\treearg{i}{t}$.
Each agent also maintains and updates a label $\labelarg{i}{t}$
indicating which problem in $\treearg{i}{t}$ it is solving.
The candidate optimal solution of a generic problem
$\probarg{i}{\ell}$ of $\treearg{i}{t}$, for some $\ell$,
is characterized in terms of a small, representative set of columns
called \emph{basis} (c.f. Appendix~\ref{sec:lps}).
We denote as $\basisarg{i}{t}$ the candidate
optimal basis of agent $i$ at time $t$.
At each
communication round $t$, the generic agent $i$ constructs
a local restricted master program \ac{RMP}$_i$ in the form
			\begin{align}
			\label{eq:rmp_alg}
	          \begin{split}
				\max_{\lambdaRMP_i} \: & \: \bar{\cc}_{V,i}^\top \lambdaRMP_i
				  \\
				\subj \: & \: \bar{V}_i\lambdaRMP_i = \1{\Ntask},
				  \\
				& \: (\1{\Nag} \1{|\lambdaRMP_i|}^\top)\lambdaRMP_i =\1{\Nag},
				  \\
				  & \: \lambdaRMP_i\geq \0{|\lambdaRMP_i|}.
			\end{split}
			\end{align}
Notice that this problem has the same structure
as~\eqref{eq:RMP} where the columns are the ones of the bases $\basisarg{j}{t}$
with $j\in \innbrs_{i,t}$.  To streamline the notation, we denote by
$\bar{V}_i$ the stack of vertexes $\extr_i$ received by the agent, by
$\bar{\cc}_{V,i}$ the stack of related costs $\cc_i^\top\extr_i$ and by
$\lambdaRMP_i$ the optimization variable.
Agent $i$ solves its local \ac{RMP}$_i$, updates the candidate basis $\basisarg{i}{t}$, and
recovers the associated optimal dual variables $[\pi_i^t \;\; \mu_i^t]$.
With the dual solution of the local \ac{RMP}$_i$ at hands, agent $i$ solves
a pricing problem
  \begin{align}\label{eq:pricingalg}
    \begin{split}
      \max_{\varz_i} \: & \: (\cc_i-\pi_i^t)^\top \varz_i
      \\
      \subj \: & \: \varz_i \in \PP_i^t,
    \end{split}
  \end{align}
  which has the same structure as~\eqref{eq:BinarySubProblem}. As discussed in
  Section~\ref{sec:CentrGAP}, this allows agents to generate
a new column $h_i$.\footnote{Here, $\PP_i^t=\{\varz_i\in\{0,1\}^\Ntask \mid \varz_i\in\PP_i, \varz_i\in\Delta_i^\ell \}$,
with $\Delta_i^\ell$ being the set of branching binary constraints associated to the problem $\probarg{i}{\ell}$ that agent $i$ is solving at iteration $t$.}
If such
column improves the overall cost, i.e., it has positive reduced cost, agent $i$
substitutes one column of $\basisarg{i}{t}$ with $h_i$.
This is done according to a so called \textsc{Pivot} operation.

Each time an agent detects convergence, or receives a label $\labelarg{j}{t}>\labelarg{i}{t}$ from some neighbor $j\in\innbrs_{i,t}$, it sets $\labelarg{i}{t+1}=\labelarg{i}{t}+1$.
Then, it retrieves the local cost and solution $\JLP,\varzLP$ from $\basisarg{i}{t+1}$ through a $\solfunc$ function.
If $\JLP\geq\jarg{i}{t}$ and $\varzLP\in\{0,1\}^{\Nag\Ntask}$, it updates the local candidate optimal cost and solution as $\jarg{i}{t+1}=\JLP$, $\zinc{i}{t+1}=\varzLP$. Otherwise, it sets $\jarg{i}{t+1}=\jarg{i}{t}$, $\zinc{i}{t+1}=\zinc{i}{t}$.
If $\JLP\geq\jarg{i}{t}$ but $\varzLP\not\in\{0,1\}^{\Nag\Ntask}$ it performs a branching operation. We denote by $\branchfunc$ the routine that updates $\treearg{i}{t}$ according to a branching on $\varzLP$. Finally, the agent starts to solve a new problem, if any, by updating, through an $\upfunc$ function, the local constraint set $\PP_i^{t+1}$.
From now on we assume that the routines $\branchfunc$ and $\upfunc$ are common to all the agents.
The whole procedure is summarized in Table~\ref{alg:distr_bnp} from the perspective of agent $i$.
\begin{algorithm}[!h]
\setstretch{1.2}
	\begin{algorithmic}[0]
	\caption{Distributed Branch-and-Price Algorithm}
	\label{alg:distr_bnp}
		\StatexIndent[0] \textbf{Initialization:} $\basisarg{i}{0} = B_{H_M}$ obtained via big-$M$, incumbent cost $\jarg{i}{0}=-\infty$

		\StatexIndent[0] \textbf{Evolution:} for all $t=1,2,\ldots$\smallskip
		\StatexIndent[0.25] Receive $\basisarg{j}{t},\labelarg{j}{t}$ from $j\in \innbrs_{i,t}$
		\StatexIndent[0.25] \textsc{Case 1:} For each $j\in\innbrs_{i,t}$, $\labelarg{j}{t}\leq\labelarg{i}{t}$
		\StatexIndent[1] Set
			\begin{align*}
	          \begin{bmatrix}
			      \bar{\cc}_{V,i}^\top
			      \\
			      \bar{V}_i
			   \end{bmatrix}\!\!
			    =%
			    \bigcup_{j\in\innbrs_{i,t} \bigcup\{i\}} \basisarg{j}{t}.
	      	\end{align*}
		\StatexIndent[1]
		  	Find optimal basis $\basisarg{i}{t+1}$ and dual solution $[\pi_i^t \;\; \mu_i^t]$ of
			\begin{align*}
	          \begin{split}
				\max_{\lambdaRMP_i} \: & \: \bar{\cc}_{V,i}^\top \lambdaRMP_i
				  \\
				\subj \: & \: \bar{V}_i\lambdaRMP_i = \1{\Ntask},
				  \\
				& \: (\1{\Nag} \1{|\lambdaRMP_i|}^\top)\lambdaRMP_i =\1{\Nag},
				  \\
				  & \: \lambdaRMP_i\geq \0{|\lambdaRMP_i|}.
			\end{split}
			\end{align*}
		\StatexIndent[1] Generate column $\colg$ solving
		\begin{align*}%
		\begin{split}
		  \max_{\varz_i} \: & \: (\cc_i-\pi_i^t)^\top \varz_i
		  \\
		  \subj \: & \: \varz_i \in \PP_i^t.
		\end{split}
		\end{align*}
		\StatexIndent[1] Update
		    $\basisarg{i}{t+1}=\textsc{Pivot} ( \basisarg{i}{t+1},\colg )$
		\StatexIndent[1]  $\jarg{i}{t+1}\!=\!\jarg{i}{t},\zinc{i}{t+1}\!=\!\zinc{i}{t},\PP_i^{t+1}\!=\!\!\PP_i^{t},\labelarg{i}{t+1}\!=\!\labelarg{i}{t},\treearg{i}{t+1}\!=\!\treearg{i}{t}$
		\StatexIndent[1] If $\basisarg{i}{t+1}$ has not changed for $2NL+1$ rounds
		\StatexIndent[1.75]
		\textsc{GOTO Case 2}
		\smallskip
		\StatexIndent[0.25] \textsc{Case 2:} There exists $j\in\innbrs_{i,t}$ s.t. $\labelarg{j}{t}>\labelarg{i}{t}$
\StatexIndent[1]
		$\labelarg{i}{t+1}=\labelarg{i}{t}+1$
		\StatexIndent[1] $\varzLP,\JLP = \solfunc(\basisarg{i}{t+1})$
    \StatexIndent[1] \textsc{Case 2.1:} $\varzLP \in \{0,1\}^{\Nag\Ntask}$, $\JLP \geq  \jarg{i}{t}$
    \StatexIndent[1.75]  $\jarg{i}{t+1}=\JLP$, $\zinc{i}{t+1}=\varzLP$
    \StatexIndent[1] \textsc{Case 2.2:} $\varzLP \not\in \{0,1\}^{\Nag\Ntask}$, $\JLP \geq \jarg{i}{t}$
    \StatexIndent[1.75]  $\treearg{i}{t+1}=\branchfunc(\treearg{i}{t},\varzLP)$
    \StatexIndent[1] If $\treearg{i}{t+1}$ is empty
		    \StatexIndent[1.75] \textsc{Halt}
  \StatexIndent[1] $\PP_i^{t+1}=\upfunc(\treearg{i}{t+1})$
	\end{algorithmic}
\end{algorithm}

The convergence properties of the Distributed Branch-and-Price algorithm are stated in the next theorem.
\begin{theorem}\label{thm:distr_bnp_conv}
Let~\eqref{eq:GAP2} be feasible and
Assumption~\ref{ass:connectivity} hold.
Consider the
sequences  $\{\jarg{i}{t},\zinc{i}{t}\}_{t\geq 0}$, $i\in\until{\Nag}$
generated by the Distributed Branch-and-Price algorithm.
Then, in a finite number $\convT\in\natural$ of communication
rounds, agents agree
on a common optimal solution $\varz^\star$ with optimal
cost value $\costval^\star$ of~\eqref{eq:GAP2},
i.e.,
$\jarg{i}{t} = \costval^\star$ and
$\zinc{i}{t} = \varz^\star$, $\forall i\in\until{\Nag}$ and $\forall t\geq \convT$.~\oprocend
\end{theorem}
We refer the reader to Appendix~\ref{sec:Analysis} for the proof of Theorem~\ref{thm:distr_bnp_conv}.
We discuss some interesting features of the proposed distributed scheme.
First, agents do not need to know the universal slotted time $t$. That is, agents
can run the steps of the distributed algorithm according to their own local clock.
If an agent is performing its computation it is assumed not to have outgoing edges
on the communication graph and the steps are performed accordingly to the available
in-neighbor bases.
This implies that the proposed distributed scheme works under
\emph{asynchronous} communication networks.
Second, as it will be shown in the analysis, the $i$-th agent can detect that convergence to an optimal basis has occurred if its
basis $\basisarg{i}{t}$ does not change for $2L\Nag+1$ communication rounds. In this way, it can halt the
steps in \textsc{Case 1} of Algorithm~\ref{alg:distr_bnp}.
Third, during the first iterations an agent $i$ may not have enough information to solve the \ac{RMP}$_i$~\eqref{eq:rmp_alg}. Thus,
it plugs into the local problem a set
of artificial variables, eventually discarded during the evolution of the algorithm,
with high cost. This method, also called Big-M method,
allows the agents to always find a solution to the \ac{RMP}$_i$.
As for the communication overhead, at each communication round each robot
  sends to its neighbors a matrix of size $(N+M+1)\times(N+M)$. Each column of
  this matrix is in the form
  $\colg = [c_i^\top \bar{v}_i \;\; \bar{v}_i^\top \;\; e_i^\top]^\top$.  Here,
  $c_i^\top \bar{v}_i$ is a real number specifying the cost to execute an
  allocation $\bar{v}_i\in\{0,1\}^M$.  The vector $e_i^\top\in\{0,1\}^N$
  specifies which robot generated that allocation. It is worth noting that the
  vector $\bar{v}_i$ can be encoded as an array of $M$ bits while $e_i^\top$ can
  be encoded as an integer number.  Finally, we underline that the assumption
that~\eqref{eq:GAP2} is feasible can be relaxed to include unfeasible GAPs, but
this assumption allows us to lighten the discussion.

\begin{remark}
As a possible variation, agents may harness the communication with a
    \emph{Cloud} node to speed-up the convergence time, and reduce the local
    memory and computing requirements.
    In this architecture the cloud unity is only involved in the storage of the
    branching tree (and not in the column generation steps). Thereby, agents do
    not construct local branching trees.  Also, agent data remain private and
    the number of messages exchanged at each communication round does not
    increase.
    When an agent $i$, at time $\bar{t}_\ell$, detects that convergence to an
    optimal solution of a problem $\Prob^\ell$ has occurred, it sends the basis
    $\basisarg{i}{\bart{\ell}}$ to the Cloud.
    At this point, the cloud extracts the optimal cost and solution $\jell$ and
    $\varzell$ and analyzes them according to the steps in \textsc{Case $2$} of
    Algorithm~\ref{alg:distr_bnp}.  Finally, if the tree is not empty, it
    extracts a new problem $\Prob^{\ell+1}$ from the tree according to the
    extraction strategy, and broadcasts to each agent $i$ the additional
    constraints to build-up $\PP_i^{\ell+1}$.
    The proof of the cloud-based version follows similar arguments as the one of
    Theorem~\ref{thm:distr_bnp_conv} and is omitted.~\oprocend
\end{remark}

\section{Numerical Computations}
\label{sec:Numerical}

In order to assess the performance and highlight the main features of our
distributed algorithm, we provide a set of numerical computations.
Simulations have been implemented on the DISROPT~\cite{farina2020disropt} toolbox and
carried out on a laptop equipped with a $2.5$ GHz dual core processor and $16$ GB of RAM.
In the following,
we generate new problems, during the branching procedure,
by adding constraints in the form $\varz_{i_k}=0$ and $\varz_{i_k}=1$, where $\varz_{i_k}$ is
the first non-integer entry of the vector $\varz$.
Regarding the order in which problems are extracted and solved,
we adopt the widely used \emph{depth first}
selection procedure,~\cite{martello1990knapsack}.
In this approach, the generated problems are
stored in a \emph{stack}, thereby the extraction procedure follows a LIFO approach.
Each time a branching occurs, the new problems are placed on the top of the stack.
In our implementation, we insert in the first position
the problem in which $\varz_{i_k}=0$ is added at last.

We perform Monte Carlo simulations
on random \ac{GAP} instances. We generate such instances according to four different random
models, usually referred to as Model A, B, C and D, of increasing difficulty. We refer the reader to
\cite{savelsbergh1997branch} for a survey on such models.
Let $\mathcal{U}(a,b)$ denote the discrete uniform distribution on the interval $[a,b]$. The data are generated as follows.
\begin{itemize}
\item \emph{Model A:} $\wel^A_{\indexa \indexb}\in\unif{10}{25}$, $\pel^A_{\indexa \indexb}\in\unif{5}{25}$ and $\vel^A_\indexa  =
9(\Ntask /\Nag) + 0.4 \max_{1\leq \indexa  \leq \Nag} \sum_{\indexb\in \mathcal{J}_\indexa^\star} \wel_{\indexa \indexb}$, with
$\mathcal{J}_\indexa^\star:=\{\indexb\mid \indexa =\argmin_{r} \pel_{r\indexb} \}$.
\item \emph{Model B:} $\wel^B_{\indexa \indexb}=\wel^A_{\indexa \indexb}$, $\pel^B_{\indexa \indexb}=\pel^A_{\indexa \indexb}$ and $\vel^B_\indexa=0.7 \vel^A_\indexa$.
\item \emph{Model C:} $\wel^C_{\indexa \indexb}=\wel^A_{\indexa \indexb}$, $\pel^C_{\indexa \indexb}=\pel^A_{\indexa \indexb}$ and $\vel^C_\indexa =
\sum_{1\leq \indexb \leq \Ntask}\wel_{\indexa \indexb}/m$.
\item \emph{Model D:} $\wel^D_{\indexa \indexb}\in\unif{1}{100}$, $\pel^D_{\indexa \indexb} =
100 - \wel_{\indexa \indexb} + k$, with $k\in\unif{1}{21}$ and $\vel^D_\indexa =\vel^C_\indexa$.
\end{itemize}

We consider different scenarios by varying the number of agents and tasks, thus considering
problems with different size and task-over-agents ratio. As for the number of agents,
$\Nag= 5,10,15$, while, for the number of tasks, $M=20,30$.

We generate $50$ random instances for each scenario and for each model.  We
are interested in both time and memory performance of the distributed algorithm.
Thus, we show the time that is needed to terminate the algorithm, expressed in
terms of the number of communication rounds, and the maximum number of tree
nodes stored by the agents. We also show the equivalent time, in seconds,
  needed for each simulation.  Since DISROPT exploits the MPI protocol to
  simulate the agents, the computation time per-agent is evaluated as
  $T_{ag} = T_{el}N_{co}/N$ where $T_{el}$ is the total elapsed time and
  $N_{co}$ is the number of cores.
As the problem size increases, the solution of these problems requires the exploration of
thousands of tree nodes, see, e.g., \cite{savelsbergh1997branch}.
However, in practical scenarios where assignment problems have to be solved almost in realtime,
it is useful to consider a feasible sub-optimal solution to the problem instead of an optimal one.
Thereby, even though our algorithm is able to find an optimal solution, in the
proposed simulations
agents interrupt the distributed algorithm when they find a feasible (sub-optimal) solution.
For this reason, we also provide the relative error, in terms of cost value, between the
exact solution (evaluated through a centralized solver) and the solution
found by the agents.
As for the connectivity among agents, we consider a static network modeled by a cyclic digraph.
We underline that our algorithm adapts to more complex graph models.
However, the choice of such digraph is interesting for simulation purposes
due to the fact that it is the
static digraph with largest diameter. Thus, the expected
number of communication rounds to completion is expected to be higher with respect to
graphs with smaller diameter.

\begin{figure}
\centering
	\includegraphics[width=.8\columnwidth]{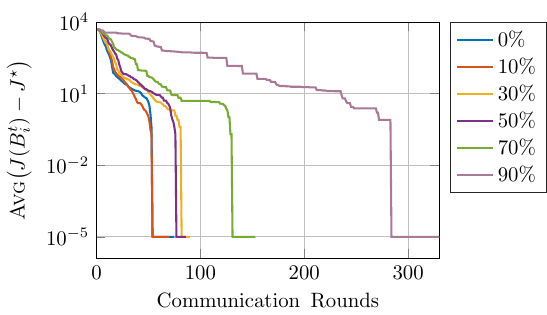}
	\caption{Cost error during the evolution of the algorithm with different percentages of packet loss.}\label{fig:lossy}
\end{figure}

The mean value and the standard deviation (evaluated over the number of
  trials) for each simulation scenario are shown in Table~\ref{tb:results}.
We highlight that, in all the simulations, the average
relative error is always below $5\%$.
The time to convergence increases with the task-to-agent ratio ($\Ntask/\Nag$). As an example, see Table~\ref{tb:results}, Model A with $\Nag=15$ and $\Ntask=30$ requires less communication rounds than Model A with $\Nag=5$ and $\Ntask=30$, even though the overall number of optimization variable is larger.
This behavior of the Distributed Branch-and-Price algorithm appears to be consistent with the one reported in the literature for centralized methods.
Similarly, Model D is far more difficult to be solved than Model A and requires
more communication rounds (see, e.g, the communication rounds needed to
  solve Model A and Model D with $\Nag=5,\Ntask=20$).  We underline that the
number of communication rounds strictly depends on the graph diameter. Since we
run the algorithm on a cyclic digraph, whose diameter is $N-1$, the results
provided in Table~\ref{tb:results} are the ones expected in case of loose
connectivity.  The maximum number of stored nodes exhibits a similar
behavior. That is, as the task-to-agent ratio increases and more difficult
models are considered, the distributed algorithm has to explore more branches.
To conclude, we propose a numerical simulation in which robots communicate in a network subject to packet loss.
We consider a scenario with $\Nag=5$ and $\Ntask=20$. Problem data are generated according to Model A.
We consider the cases with loss probability $0\%$ (no packet loss), $10\%,30\%,50\%,70\%,90\%$. Specifically, at each iteration, the $i$-th robot discards the message from the $j$-th robot according to the given probability.
 Results are given in Figure~\ref{fig:lossy}.  We show the mean error between the cost $J(\basisarg{i}{t})$ associated to the basis $\basisarg{i}{t}$ and the optimal solution $J^\star$. 

\begin{table*}[ht]
\footnotesize
	\begin{center}
	\caption{Numerical Results}\label{tb:results}
		\begin{tabular}{ccccccc}%
			Model & N & M & Communication Rounds  (Avg--Std) & Relative Error (Avg--Std) & Stored Nodes (Avg--Std) & Time (Avg--Std)\\ \hline
			\multirow{6}{*}{A}%
&$5$ & $20$ & $83.30$--$25.36$ & $0.00 \%$--$0.00 \%$ & $1.10$--$0.36$ & $1.94$--$0.60$\\
&$5$ & $30$ & $329.38$--$151.87$ & $0.01 \%$--$0.08 \%$ & $1.44$--$0.64$ & $12.85$--$5.89$\\
&$10$ & $20$ & $75.34$--$43.86$ & $0.00 \%$--$0.00 \%$ & $1.30$--$0.78$ & $3.34$--$1.86$\\
&$10$ & $30$ & $107.92$--$87.52$ & $0.01 \%$--$0.07 \%$ & $1.30$--$1.06$ & $10.29$--$7.88$\\
&$15$ & $20$ & $76.60$--$26.23$ & $0.01 \%$--$0.05 \%$ & $1.12$--$0.38$ & $5.09$--$2.01$\\
&$15$ & $30$ & $95.86$--$33.27$ & $0.00 \%$--$0.00 \%$ & $1.08$--$0.34$ & $14.59$--$4.80$\\
			\hline
			\multirow{6}{*}{B}%
& $5$ & $20$ & $192.04$--$205.30$ & $1.06 \%$--$2.10 \%$ & $3.50$--$3.79$ & $4.67$--$5.03$\\
& $5$ & $30$ & $774.50$--$680.58$ & $0.57 \%$--$0.91 \%$ & $5.04$--$4.40$ & $33.68$--$29.46$\\
& $10$ & $20$ & $161.36$--$222.29$ & $0.25 \%$--$0.80 \%$ & $3.14$--$4.38$ & $6.56$--$9.13$\\
& $10$ & $30$ & $236.36$--$324.55$ & $0.20 \%$--$0.53 \%$ & $3.24$--$4.53$ & $23.42$--$31.24$\\
& $15$ & $20$ & $90.02$--$42.64$ & $0.02 \%$--$0.15 \%$ & $1.36$--$0.66$ & $5.16$--$2.43$\\
& $15$ & $30$ & $178.40$--$174.70$ & $0.04 \%$--$0.10 \%$ & $2.16$--$2.13$ & $27.32$--$26.81$\\
			\hline
			\multirow{6}{*}{C}%
& $5$ & $20$ & $158.24$--$149.94$ & $0.63 \%$--$1.25 \%$ & $3.00$--$3.03$ & $3.95$--$3.74$\\
& $5$ & $30$ & $652.32$--$741.28$ & $0.48 \%$--$1.20 \%$ & $4.48$--$5.59$ & $27.98$--$31.95$\\
& $10$ & $20$ & $155.06$--$230.47$ & $0.47 \%$--$1.58 \%$ & $3.24$--$4.93$ & $5.81$--$8.44$\\
& $10$ & $30$ & $375.52$--$431.19$ & $0.59 \%$--$1.03 \%$ & $5.50$--$6.41$ & $35.37$--$40.32$\\
& $15$ & $20$ & $107.08$--$133.53$ & $0.14 \%$--$0.57 \%$ & $1.80$--$2.26$ & $4.85$--$6.16$\\
& $15$ & $30$ & $294.02$--$306.51$ & $0.24 \%$--$0.42 \%$ & $3.94$--$4.15$ & $38.92$--$39.99$\\

			\hline
			\multirow{6}{*}{D}%
			& $5$ & $20$ & $1072.76$--$601.48$ & $4.31 \%$--$3.88 \%$ & $20.10$--$11.23$ & $33.47$--$18.10$\\
& $5$ & $30$ & $4805.44$--$2525.97$ & $4.91 \%$--$3.65 \%$ & $36.80$--$18.05$ & $201.56$--$101.85$\\
& $10$ & $20$ & $933.76$--$979.26$ & $2.77 \%$--$4.12 \%$ & $18.44$--$20.06$ & $36.20$--$36.08$\\
			&$10$ & $30$ & $5959.95$--$3596.71$  & $4.96\%$--$4.26\%$ &  $63.55$--$38.69$ & $560.14$--$338.02$\\
			& $15$ & $20$ & $187.88$--$164.20$ & $0.22 \%$--$0.44 \%$ & $2.84$--$2.52$ & $9.38$--$8.19$\\
			&$15$ & $30$ & $6171.65$--$4015$ & $3.37\%$--$2.81\%$ & $56.15$--$37.02$ & $740.5$--$481.8$\\
			\hline
		\end{tabular}
	\end{center}
\end{table*}

\begin{remark}
  Other distributed approaches suitable for the GAP solution are the ones
    in~\cite{choi2009consensus,luo2013distributed}.  As for the one
    in~\cite{choi2009consensus}, authors consider the case in which
    $\wel_{i\indexb} \in \{0,1\}$ for each $i\in\until{\Nag}$ and for each
    $\indexb\in\until{\Ntask}$ (cf.~\eqref{eq:GAP}). When the cost function is
    linear, as in the GAP scenario, the constraint matrix is said to be totally
    unimodular and the problem can be solved as a linear problem instead of a
    mixed-integer problem. Thereby, the first solution found by our algorithm, which is also tailored for general GAPs with non-unimodular structure, is
    always the optimal one. The one found by the scheme
    in~\cite{choi2009consensus} is guaranteed to be at most $50\%$ suboptimal.
    Moreover, our algorithm allows for directed communication graphs, while the
    one in~\cite{choi2009consensus} assumes undirected communications.  Finally,
    agents in~\cite{choi2009consensus} exchange, at each communication round,
    two real vectors of size $\Nag$ and $\Ntask$ respectively and a vector of
    size $\Ntask$ representing which agent is performing each task.  As for the
    distributed approach in~\cite{luo2013distributed}, each agent has to flood
    its local variables to all the other agents.  This results in multi-hop
    communications at each iteration. Moreover, the scheme
    in~\cite{luo2013distributed} is based on the assumption of static,
    undirected graphs.  Each agent in~\cite{luo2013distributed} sends to its
    neighbors a real vector of size $M$ and three integers. Similarly to our
    approach, it considers the solution of a knapsack problem at each iteration.
    Finally, as~\cite{choi2009consensus}, it guarantees at most $50\%$
    sub-optimality of the solution found.
We perform a comparison between the proposed approach and the one in~\cite{luo2013distributed} for the scenario with $\Nag=15,\Ntask=30$. The results are in Table~\ref{tb:comparison}. We took for both the schemes  the same underlying communication graph. Since the algorithm in~\cite{luo2013distributed} needs agents to flood their information to all the other agents at each communication round, we multiply the total number of iterations of the algorithm by $\Nag d$, with $d$ diameter of the graph.
Besides the problems generated via Model D, our algorithm is able to find in less iterations a solution with a smaller relative error with respect to the one found by the algorithm in~\cite{luo2013distributed}.~\oprocend 

\end{remark}

\begin{table*}[h]\caption{Performance Comparison}\label{tb:comparison}
\footnotesize
\begin{center}
\begin{tabular}{ccc|cc|cc}
\hline
\multirow{2}{*}{Model}&\multirow{2}{*}{M}&\multirow{2}{*}{N}&\multicolumn{2}{c|}{Distributed Branch-and-Price} & \multicolumn{2}{c}{[33]} \\
\cline{4-7}
& &  & {Comm. Rounds}  (Avg--Std) & Rel. Error (Avg--Std) & Comm. Rounds  (Avg--Std) & Rel. Error (Avg--Std) \\
\hline
			\multirow{1}{*}{A}%
&$15$ & $30$ & $95.86$--$33.27$ & $0.00 \%$--$0.00 \%$ &  {$411.6$--$62.44$} &  {$0.56\%$--$0.37\%$}\\
			\hline
			\multirow{1}{*}{B}%
& $15$ & $30$ & $178.40$--$174.70$ & $0.04 \%$--$0.10 \%$ &  {$447.3$--$65.77$} &  {$1.91\%$--$1.13\%$}\\
			\hline
			\multirow{1}{*}{C}%
& $15$ & $30$ & $294.02$--$306.51$ & $0.24 \%$--$0.42 \%$ &  {$453.6$--$107.9$} &  {$3.45\%$--$1.74\%$}\\

\hline
\multirow{1}{*}{D}%
&$15$ & $30$ & $6171.65$--$4015$ & $3.37\%$--$2.81\%$ & $302.4$--$34.12$ & $0.0\%$--$0.0\%$\\
\hline
\end{tabular}
\end{center}
\end{table*}

\section{Experiments on GAPs\\ for a Team of Ground and Aerial Robots}
\label{sec:experiment}
In the following, we provide experimental results on a generalized
  assignment scenario where a team of heterogeneous (ground and aerial) mobile
  robots has to accomplish a set of tasks that may not be completely known in
  advance.
We start by describing how we implemented the proposed distributed
scheme into the ROS framework. Then, we propose the \emph{\DASS}, 
a resolution methodology for this assignment scenario,
and provide experiments on a real fleet of ground and aerial robots.

\subsection{Experimental ROS Architecture}
\begin{figure}
	\centering
	\includegraphics[width=\columnwidth]{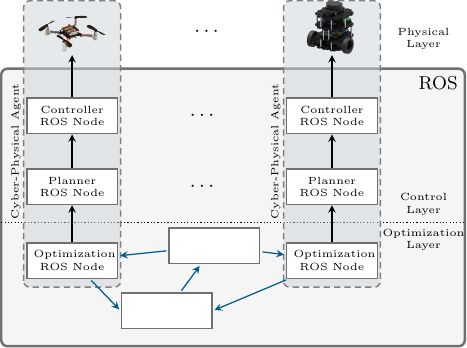}
	\caption{{\DASS} architecture. Blue rectangles represent the
	smart cyber-physical agents endowed with computation, communication and actuation capabilities.}
\label{fig:arch}
\end{figure}

In the proposed architecture, robots are ``smart'' cyber-physical agents
  endowed with communication, computation and actuation capabilities.  Each
cyber-physical agent consists of three \emph{ROS nodes}, namely
\emph{Optimization}, \emph{Planner} and \emph{Controller} ROS nodes, see
Figure~\ref{fig:arch}.
It is worth noticing that, in general, each agent has a dedicated machine on which these processes run, so that there is no need for
a central computing unit handling the agents.
The Optimization node handles the steps of the distributed optimization algorithm of the associated cyber-physical agent.
It communicates with the Optimization nodes of the other robots through the ROS \emph{publisher-subscriber} communication protocol
according to a fixed communication graph, and exchanges messages containing the local candidate bases.
Note that the communication among processes in ROS is completely \emph{asynchronous}.
As shown in the theory this is handled by our distributed algorithm.
Each time such process receives a message from a neighbor,
a \emph{callback function} stores the received
basis. Each node performs an iteration of the Distributed Branch-and-Price algorithm
within a \emph{loop} of $5$ ms. At the beginning of this loop,
the node performs one step of the column generation
algorithm with the received bases.
Then, it sends the updated basis to its neighbors and stays idle
until the next loop iteration.
The Optimization nodes
characterize the \emph{Optimization Layer} (c.f. Figure~\ref{fig:arch}) of the proposed
architecture.
The Control and Planner ROS nodes constitute instead the \emph{Control Layer} of the
proposed software.
More in detail, the Planner node generates,
through polynomial splines, a sufficiently smooth trajectory
steering a robot over its designated tasks.
The Controller implements a trajectory tracking strategy. It receives the pose
of the vehicle by a Vicon motion capture system and sends the control inputs to
the robot actuators (\emph{Physical Layer} in Figure~\ref{fig:arch}).

\subsection{Distributed Dynamic Assignemnt: Scenario and Strategy}
\begin{figure*}[ht]
	\centering
	\includegraphics[width=.7\textwidth]{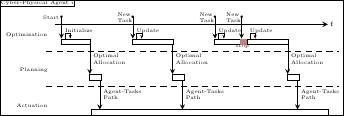}
	\caption{An example of the {\DASS} evolution from the perspective of the generic cyber-physical agent. Each time a new task appears, the robot updates the local problem data and re-starts the optimization. If a new task arrives during the re-optimization, this latter is halted (red rectangle) and a new one starts. When robot-to-task paths are evaluated, robot actuation changes accordingly.}
	\label{fig:dynamic_flow}
\end{figure*}

\begin{figure}[t]
    \centering
\includegraphics[width=.98\columnwidth]{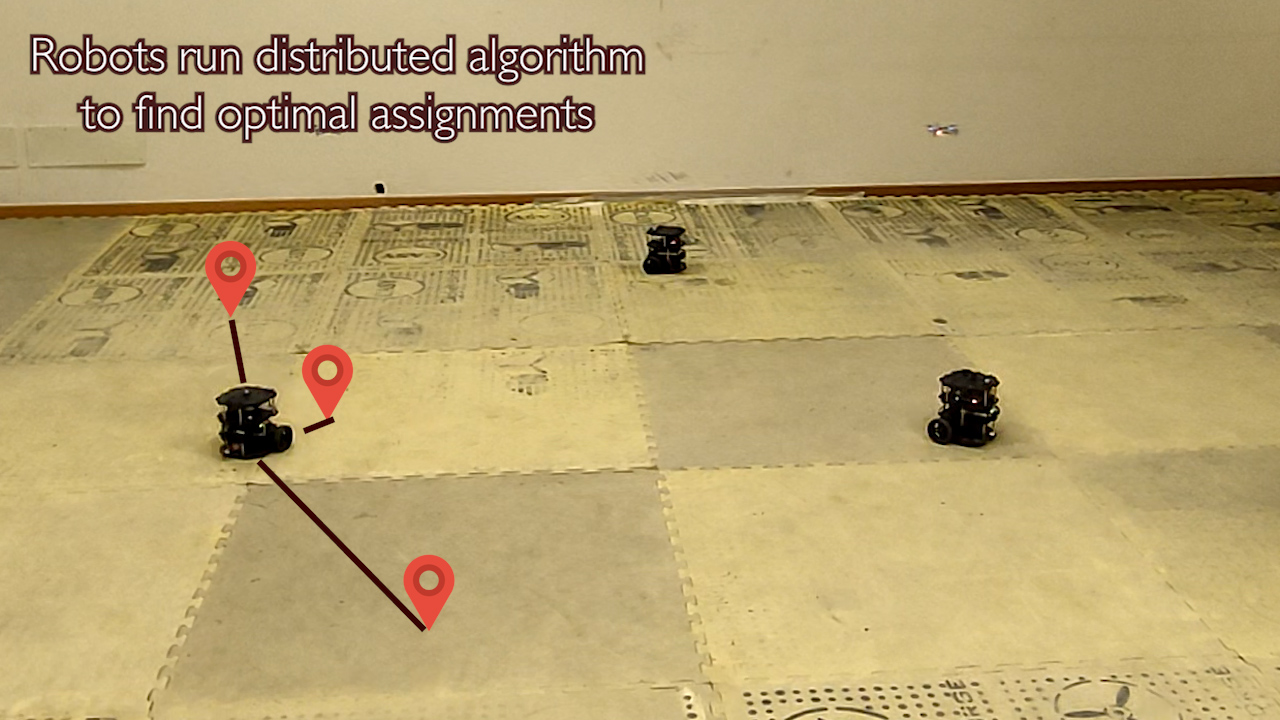}
    \caption{Snapshot from an experiment. The figure depicts the optimal assignment for one of the robots.}
    \label{fig:experiment_gap}
\end{figure}

The scenario evolves as follows.  We consider a team of ground and aerial mobile
robots moving in a three-dimensional environment parametrized by a frame
$\{x,y,z\}$.
A set of tasks, parametrized by a position on the $\{x,y\}$ plane, are scattered in the environment.
Some of the tasks can be accomplished only by ground robots, other are accessible only to aerial robots 
and there are tasks that can be performed by all the robots.
For a task to be accomplished, a robot has to visit the task location, stand
  still for a certain random time $T^{\textsc{h}}$ and go back to a given depot
  (e.g., to recharge batteries).  As in practical applications, the information
about the problem instance is not known in advance and new data arrive while the
agents are fulfilling other tasks.
To adapt the Distributed Branch-and-Price algorithm to such dynamic scenario, we
combine it with the methodology proposed in~\cite{karaman2008large} into an
optimization and task-fulfilling approach which we call Distributed Dynamic
Assignment and Servicing Strategy.  Such procedure combines a distributed
optimization phase with a planning and control scheme to steer the robots over
the assigned tasks.  More in detail, the experiment starts with the
cyber-physical agents running
the Distributed
Branch-and-Price Algorithm on a set of tasks known in advance.
Inspired by~\cite{choi2009consensus}, we pick
$\pel_{im}$ in~\eqref{eq:GAP} as a \emph{ time-discounted reward},
i.e., $\pel_{im}=\lambda_m^{\tau_i^m}$ where
$\lambda_m\in(0,1)$ is a scoring value for task $m$ and $\tau_i^m$ is the time needed
by agent $i$ to reach task $m$.
The time $\tau_i^m$ is evaluated as the robot-task distance (on the $\{x,y\}$ plane)
scaled by the robot maximum speed ($1 \textrm{ m/s}$ for the UAVs and $0.22 \textrm{ m/s}$ for the ground vehicles).
The fact that a task $m$ is not accessible to a certain robot $i$ is modeled by taking $\wel_{im}>g_i$ in~\eqref{eq:GAP}.
In the following, we assume that the sets $\PP_i$, generated randomly according to Model A
in Section~\ref{sec:Numerical},
are fixed throughout the scenario evolution.
As soon as a robot reaches the designed task, it stands still on the location for a random time $T^{\textsc{{h}}}$ between $3$ and $5$ seconds.
In the proposed experiment, we consider a dynamic scenario in which the number of tasks appearing during the evolution is always smaller than the number of served tasks.
For the sake of simplicity, we suppose that one new task is made available
to robots each time a task has been fulfilled.
In this way, the size of the optimization problem is constant.  We point out
that the strategy can be applied to more general cases where more tasks are
revealed. Moreover, while in the current set-up we consider the immediate
  strategy in which we re-optimize the entire problem, one could think of
  implementing tailored schemes leveraging the dynamic structure of the
  problem.  As soon as new tasks appear, the cyber-physical agents run the
Distributed Branch-and-Price algorithm on a problem including the new tasks and
discarding the visited ones.  Specifically, the cost vector entries change
according to the new task positions.  Meanwhile, each robot keeps
performing tasks according to its latest allocation.  An example of the
  evolution of this strategy is in Figure~\ref{fig:dynamic_flow}. 
  A snapshot
  from an experiment with $2$ Crazyflie nano-quadrotors and $3$ Turtlebot3
  Burger is in Figure~\ref{fig:experiment_gap}. Here robots have terminated the
  distributed optimization procedure and one of the allocations is shown.
 A video is available as supplementary material to the paper\footnote{The video is also available at \texttt{\url{https://youtu.be/Sl_3ZmJvvbU}}.}.
\begin{remark}
  As discussed, e.g., in~\cite{fisher1981generalized}, GAPs can be also used to
  find approximate solutions of vehicle routing problems (VRPs).  In general,
  VRPs penalize the order of execution of the tasks, and involve a larger number
  of variables with respect to GAPs. The idea in~\cite{fisher1981generalized} is
  to construct a GAP instance based on the VRP problem data. As soon as a GAP
  solution has been found, robots perform their associated task in an order that
  minimizes, e.g., the total travelled distance. This can be done, e.g., by
  solving a Shortest Hamiltonian Path Problem (SHPP).  The proposed {\DASS}
  could be thus modified in order to address such scenarios. Specifically,
  robots start solving the GAP with the available tasks and, once an optimal
  solution has been found, construct robot-to-tasks paths by solving SHPPs. When
  a new task arrives, robot re-solve the optimization problem and adjust the
  path according to the new problem data.  We performed an experiment with $3$
  Crazyflie nano-quadrotors and $2$ Turtlebot3 Burger with the cloud-based
  approach.
  A video is available as supplementary material to the manuscript.\footnote{The video is also available at~\texttt{\url{https://youtu.be/vBSJsduFYKQ}}.}~\oprocend
\end{remark}

\section{Conclusion}
In this paper, we proposed a purely distributed branch-and-price approach to solve the Generalized Assignment Problem
in a network of agents, endowed with computation and communication capabilities, that are aware of only a small part
of the global optimization problem data.
Agents cooperatively solve a relaxations of the GAP by means of a distributed column generation algorithm, targeted for this
particular scenario involving binary optimization variables. Since the solution of this relaxation may not be feasible for the GAP,
agents cooperatively generate and solve new optimization problems, considering each time additional constraints.
Finally, we considered an assignment scenario where tasks may appear dynamically during time. We implemented
the proposed algorithms in a ROS based testbed and showed results from experiments on a team of ground and aerial vehicles executing
the generalized assignment.
Future investigations may include the solution of dynamic instances of the GAP with tailored distributed approaches that do not need to re-optimize the entire problem when new data arrive.

\section*{Acknowledgment}
The authors would like to thank Alessandro Rucco for the fruitful discussions and Nicola Mimmo for the support during the experiments.

\appendices
\section{Linear Programs}\label{sec:lps}
An LP in \emph{standard form} is a problem in the form
\begin{align}\label{eq:lp_primal}
\begin{split}
\min_{x} \: & \: c^\top x
  \\
\subj \: & \: Ax=b,
  \\
& \: x\geq 0.
\end{split}
\end{align}
where $c\in\real^d$, $A\in\real^{r\times d}$ and $b\in\real^r$
are the problem data and $x\in\real^d$ is the optimization variable.
All the problem constraints are expressed as equality constraints
and the variables must be non-negative.
A \emph{column} for the problem in~\eqref{eq:lp_primal} is a vector in the form $[c_\ell\;\; A^\top_\ell]\in\real^{r+1}$ where $A^\top_\ell$ is
the $\ell$-th column of $A$. A \emph{basis} $B$ is a set of $r$ independent columns of the LP. We denote by $c_B$ ($A_B$)
the sub-vector (sub-matrix) of $c$ ($A$) constructed from the columns in $B$.
Assume that a solution $x^\star$ to~\eqref{eq:lp_primal} exists.
Then, it can be shown that $x^\star$ can be decomposed
into two sub-vectors $x^\star_B\neq 0$ of \emph{basic} variables and $x^\star_N=0$ of \emph{non-basic} variables.
A basis represents a \emph{minimal representation} of a linear program, i.e.,
it is a subset of the problem data representing the problem solution.
It can be shown that there exists a basis $B$ such that $x^\star_B$ is the solution of:
\begin{align*}
\begin{split}
\min_{x} \: & \: c_B^\top x
  \\
\subj \: & \:A_B x=b,
  \\
& \: x\geq 0.
\end{split}
\end{align*}

\section{Proof of Theorem~\ref{thm:distr_bnp_conv}}\label{sec:Analysis}
\subsection{Preliminary Lemmas for the Proof of Theorem~\ref{thm:distr_bnp_conv}}
Before proceeding with the proof of Theorem~\ref{thm:distr_bnp_conv}, we
provide two lemmas which are useful for the analysis.
\begin{lemma}\label{lemma:lp_convergence}
  Let Assumption~\ref{ass:connectivity} hold. Consider a network of $N$ agents
  running the steps in \textsc{Case 1} of Algorithm~\ref{alg:distr_bnp} to solve
  a node $\Prob^\ell$ of the tree. Then, in a finite number of iterations,
  agents reach consensus to an optimal basis $B^\ell$ associated with the
  optimal cost $J^{\star\ell}$ of $\Prob^\ell$.
\end{lemma}
\begin{IEEEproof}
  The proof mimics the one proposed in~\cite{burger2011locally}. We refer the
  reader to this work for additional details.  First, we show that
    $\Prob^\ell$ can be obtained by applying the Dantzig-Wolfe decomposition to
    the following Linear Program:
\begin{align}\label{eq:Conv-GAP}
\begin{split}
  \max_{\varz_1,\ldots,\varz_\Nag} \: & \: \sum_{i =1}^\Nag c_i^\top \varz_i
  \\
  \subj \: & \: \sum_{i=1}^\Nag \varz_i = \1{\Ntask},
  \\
  & \: \varz_i \in \text{conv}(\PP_i^\ell),  i= 1,\ldots, \Nag.
  \\
\end{split}
\end{align}
Indeed, we recall that, for GAPs, the vertexes of $\text{conv}(\PP_i^\ell)$ coincide with
the points $\extr_i, q\in\QQ_i^\ell$,~\cite{barnhart1998branch}. Thus, points $z_i\in\text{conv}(\PP_i^\ell)$
can be represented as $\varz_i = \sum_{q=1}^{|\QQ_i^\ell|} \extr_i \mult_i$ with $\sum_{q=1}^{|\QQ_i^\ell|} \mult_i = 1$ and $\mult_i\geq 0$. Let $\Lambda$ be the stack of the variables $\mult_i$. By substituting these equations in~\eqref{eq:Conv-GAP} one obtains
a problem in the form
\begin{align*}
\begin{split}
\max_{\Lambda} \: & \: \sum_{i =1}^\Nag \sum_{q=1}^{|\QQ_i^\ell|} (\cc_i^\top \extr_i)\mult_i
  \\
\subj \: & \: \sum_{i=1}^\Nag \sum_{q=1}^{|\QQ_i^\ell|} \extr_i \mult_i =\1{\Ntask},
  \\
  & \: \sum_{q=1}^{|\QQ_i^\ell|} \mult_i = 1,  i= 1,\ldots, \Nag,
  \\
  & \: \mult_i\geq 0, q\in\until{|\QQ_i^\ell|}, i= 1,\ldots, \Nag,
\end{split}
\end{align*}
which is problem $\Prob^\ell$.
Notice that the resulting pricing problem for each agent $i$ is
\begin{align}\label{eq:pricing_conv}
\begin{split}
  \max \: & \: (c_i-\pi)^\top \varz_i
  \\
  \subj \: & \: \varz_i \in \text{conv}(\PP_i^\ell).
\end{split}
\end{align}
By definition of convex hull and linearity of the cost function,~\eqref{eq:pricing_conv} shares the same optimal vertexes of~\eqref{eq:BinarySubProblem}.
Thereby, the steps of \textsc{Case 1} in Algorithm~\ref{alg:distr_bnp}
can be seen as applied to LP~\eqref{eq:Conv-GAP}.
At this point, we note that, during the algorithmic evolution,
agent $i$ can update 
its local candidate basis
$\basisarg{i}{t}$ by considering new columns in the local linear
program.
It is worth noting that, starting from any basis $B_i^t$, 
there exists a finite number of pivoting operations to the 
optimal basis $B^\ell$. 
These columns can be found in two ways:
(i) by the local column generation routine and
the subsequent pivoting
and (ii) when collecting all the in-neighbors 
matrices $\basisarg{j}{t}$ with $j\in\innbrs_{i,t}$.
If $J_i^t < J^{\star\ell}$, there always exists an agent $j$ able to 
generate a column improving the cost $J_i^t$ after a pivoting.
Since the network is connected
if that column is fundamental for the evolution of the algorithm, e.g., 
it belongs to the optimal basis, then agent $j$ will
 generate it (and include it) in its basis within a
finite number of communication rounds. Thus, as soon as $J_i^t < J^{\star\ell}$,
there always exists a finite time $T_D$ such that $J_i^t < J_I^{t+T_D}$.
Since there exists only a finite number of columns, in a finite number of
communication rounds $T_f$ it stands that $J_i^{T_f}=J^{\star\ell}$ for each $i$. If a lexicographic solver is considered, then it stands that $B_i^t=B^\ell$ for each $i$. This concludes the proof.
\end{IEEEproof}

\begin{lemma}\label{lemma:halting}
Let Assumption~\ref{ass:connectivity} hold. Then, 
a processor has computed its final basis 
and can halt the execution of the steps in \textsc{Case 1} of Algorithm~\ref{alg:distr_bnp} as soon as the value of 
$B_i^t$ has not changed after $L(2N-1)$ communication rounds.
\end{lemma}
\begin{IEEEproof}
Assume that a certain node $i$ satisfies $\basisarg{i}{t} = B^\star, \jarg{i}{t}=J^\star$ for all $t \in \{t_0, \ldots, t_0 + 2L(N-1)\}$,
and pick any other node $j$. Without loss of generality, consider $t_0 = 0$.
By $L$-strong connectivity, after at most $L$ communication rounds, agent $i$
has been able to spread its basis at least to another agent.
We now define the set
$\bar{\mathcal{N}}_0$ of agents $k\in\{1,\dots N\}$ such that there exists
an increasing sequence of time instants $\tau_0,\ldots,\tau_m$ comprised between
$0$ and $L$ (i.e., with $0 \leq \tau_0$ and $\tau_m \leq L$), such that
the edges $(i,\ell_1),\dots,(\ell_{m},k)$ belong to the digraph at times $\tau_0,\ldots,\tau_m$.
This set is not empty, since the union graph is strongly connected in $[0,L]$.
Then it stands $\jarg{k}{L} \geq \jarg{i}{t}$, $\forall k\in\bar{\mathcal{N}}_0$.
Consider now the interval $[L, 2L]$, for which we define a set similar to $\bar{\mathcal{N}}_0$,
but with paths originating from the agents in $\bar{\mathcal{N}}_0 \cup \{i\}$.
Formally, consider the set $\bar{\mathcal{N}}_1$ of agents $k\in\{1,\dots N\}$ such that there exists
an increasing sequence of time instants $\tau_0,\ldots,\tau_m$ comprised between
$L$ and $2L$ (i.e., with $L \leq \tau_0$ and $\tau_m \leq 2L$), such that
the edges $(h,\ell_1),\dots,(\ell_{m},k)$ belong to the digraph at times $\tau_0,\ldots,\tau_m$,
for some $h \in \bar{\mathcal{N}}_0 \cup \{i\}$.
Notice that $\bar{\mathcal{N}}_0\subset \bar{\mathcal{N}}_1$, so that
$\bar{\mathcal{N}}_1$ has a larger cardinality than $\bar{\mathcal{N}}_0$.
Otherwise, the graph would not be strongly connected in $[L, 2L]$.
Then it stands
$\jarg{k}{2L} \geq \jarg{i}{t}$, $\forall k\in\bar{\mathcal{N}}_1$.
Iterating at most $N-1$ times, we see that the sets $\bar{\mathcal{N}}_0 \ldots, \bar{\mathcal{N}}_{N-2}$
become larger and larger, so that $j \in \bar{\mathcal{N}}_{N-2}$.
Thus, it stands that $\jarg{k}{(N-1)L}\geq \jarg{i}{t}$.
That is, after $(N-1)L$ communication rounds, all the agents have at least the
same cost agent $i$ had at time $0$.
By repeating the same arguments for the converse path, we conclude that
$\jarg{i}{2(N-1)L}\geq\jarg{k}{(N-1)L}$.
But, by assumption, $\jarg{i}{2(N-1)L} = J^\star$, so that we conclude
$J^\star \leq \jarg{j}{(N-1)L} \leq J^\star$,
i.e., $\jarg{j}{(N-1)L} = J^\star$.
Thus, if $\basisarg{i}{t}$ does not change for $L(2N-1)$ time instants,
then its value will never change afterwards because all bases $B_{j}^t, j \in \until{N}$,
have cost equal to $J^\star$ at least as early as time equals $LN$.
\end{IEEEproof}

\subsection{Proof of Theorem~\ref{thm:distr_bnp_conv}}

In order to prove the statement, we show that there exists a monotonically increasing time sequence $\{\bar{t}_\ell\}_{\ell\in\{0,\ldots, \ell_{\text{end}} \}}$, for some $\ell_{\text{end}}\in\natural$, such that, at each $\bar{t}_\ell$:
\begin{itemize}
	\item[i)]For all $i,j\in\until{\Nag}$, $\treearg{i}{\bar{t}_\ell}=\treearg{j}{\bar{t}_\ell}$, $\probarg{i}{\ell}=\probarg{j}{\ell}=\Prob^\ell$,  $\labelarg{i}{\bart{\ell}}=\labelarg{j}{\bart{\ell}}=\ell$ and there exists some $i\in\until{\Nag}$ such that $\labelarg{i}{\bart{\ell}-1}\neq\ell$;
	\item[ii)] in a finite number of communication rounds, at a time $\bart{\ell+1}\leq\bart{\ell}+Q^\ell$, $Q^\ell\in\natural$, either $\probarg{i}{\ell+1}= \probarg{j}{\ell+1}$ (with $\treearg{i}{\bar{t}_{\ell+1}}=\treearg{j}{\bar{t}_{\ell+1}}$ and $\labelarg{i}{\bart{\ell+1}}=\labelarg{j}{\bart{\ell+1}}, \, \forall i,j\in\until{\Nag}$) or agents halt the distributed algorithm, i.e., $\ell=\ell_{\text{end}}$, with $\jarg{i}{\bar{t}_{\ell+1}}=J^\star$ and $\zinc{i}{\bar{t}_{\ell+1}}=z^\star$ optimal cost and solution of~\eqref{eq:GAP2} for all $i\in\until{\Nag}$.
\end{itemize}
First notice that i) holds trivially at $t_0=0$, since all the agents start solving the relaxed version of~\eqref{eq:GAP_IPMP}, namely $\Prob^0$, and each agent initializes $\labelarg{i}{0}=0$.
Now, we assume that i) holds for some $\ell$ and prove that ii) holds.
Then, by applying the arguments
in Lemma~\ref{lemma:lp_convergence} %
agents reach consensus, in a finite number of communication rounds $\bar{Q}^\ell$,
on a basis $B^\ell$ corresponding
to an optimal solution $\lambdaell$ of $\Prob^\ell$.
Moreover, by Lemma~\ref{lemma:halting}
each agent $i$ can halt, at some time $\bar{Q}^\ell\leq t_{i,\ell}\leq \bart{\ell} + \bar{Q}^\ell+2L\Nag+1$,
the steps of
\textsc{Case 1} if its basis $\basisarg{i}{t}$ has not changed for $2L\Nag+1$ communication rounds (c.f. Algorithm~\ref{alg:distr_bnp}).
At these times, each agent obtains the same cost $\jell$ and solution $\varzell$ of $\Prob^\ell$ (retrieved from $\lambdaell$ by applying~\eqref{eq:comb_zi}) and sets
$\labelarg{i}{t_{i,\ell}+1}=\ell+1$.
If \textsc{Case 2.1} in Algorithm~\ref{alg:distr_bnp} occurs,
then each agent $i\in\until{\Nag}$ sets
$\jarg{i}{t_{i,\ell}+1}=\jell$
and $\zinc{i}{t_{i,\ell}+1}=\varzell$.
Instead, if \textsc{Case 2.2} occurs, each agent
expands the local tree $\treearg{i}{t_{i,\ell}+1}$.
Notice that agents run the $\branchfunc$ routine on the same data ($\treearg{i}{\bart{\ell}}$ and $\varzell$), so they update the same tree with the same new problems.
Finally, if there are still problems to be solved in $\treearg{i}{t_{i,\ell}+1}$, each agent extracts a new problem $\probarg{i}{\ell+1}$.
Since the routine $\upfunc$ is common to all
the agents and the constructed
trees are identical, $\probarg{i}{\ell+1}=\Prob^{\ell+1}$ for all $i$.
Otherwise, if $\treearg{i}{t_{i,\ell}+1}$ is empty, each agent halts the
Distributed Branch-and-Price Algorithm.
Let $Q^\ell=\bar{Q}^\ell+2L\Nag+2$ and let
$\bart{\ell+1}=\max_i\{t_{i,\ell}\}+1$.
From the above arguments,
$\bart{\ell}\leq\bart{\ell+1}\leq\bart{\ell}+Q^\ell$.

Now we show that, if agents halt the distributed algorithm, i.e., $\ell=\ell_{\text{end}}$, then
$\jarg{i}{\bart{\ell+1}}=J^\star$ and $\zinc{i}{\bart{\ell+1}}=\varz^\star$  for all $i\in\until{\Nag}$,
with $J^\star$ and $z^\star$ optimal cost and solution of~\eqref{eq:GAP2}.
First, notice that $\jarg{i}{\bart{\ell+1}}$ and $\zinc{i}{\bart{\ell+1}}$ are the optimal cost value and solution of some problem $\Prob^\ell$ such that $\jell\geq\jarg{i}{t}$ for each $t\leq \bart{\ell+1}$ and $\varzell\in\{0,1\}^{\Nag\Ntask}$.
Since~\eqref{eq:GAP2} is feasible, and in the branch-and-price algorithm all the nodes of the tree are explored (except the ones discarded during the pruning operation) then each agent has run at least one time the steps in \textsc{Case 2.1}.
Thereby, $\jarg{i}{\bart{\ell+1}}$ must be equal to the optimal
cost value $J^\star$ of~\eqref{eq:GAP2} and, similarly,
$\zinc{i}{\bart{\ell+1}}=\varz^\star$ optimal solution to~\eqref{eq:GAP2}.
To conclude, we underline that agents can generate only a finite number of problems.
Indeed, the number of additional constraints ($\varz_{i_k}=0$ and $\varz_{i_k}=1$) they can add is at most $2^{\Nag\Ntask}$.
Thus, there exists a time $\convT=\bar{t}^{\ell_{\text{end}}}+Q^{\ell_{\text{end}}}$ in which all the agents must halt the
distributed scheme. This concludes the proof.
\bibliographystyle{IEEEtran}
\bibliography{bibliography}
\end{document}